\documentclass{amsart}
\usepackage{amscd,amssymb,graphicx}
\input{xy}
\xyoption{all}

\author[DeCleene]{Chris DeCleene}
\address{
Chris DeCleene\\
University of Wisconsin-Eau Claire\\
Eau Claire, WI 54702-4004}
\email{decleecj@uwec.edu}

\author[Penkava]{Michael Penkava}
\address{
Michael Penkava\\
University of Wisconsin-Eau Claire\\
Eau Claire, WI 54702-4004} \email{penkavmr@uwec.edu}

\author[Phillipson]{Mitch Phillipson}
\address{
Mitch Phillipson\\
University of Wisconsin-Eau Claire\\
Eau Claire, WI 54702-4004}
\email{phillima@uwec.edu}

\thanks{This research was supported by National Science Foundation Grant 144-2974/92974, and grants from the University of Wisconsin-Eau Claire}%

\newtheorem{thm}{Theorem}[section]

\theoremstyle{definition}
\newtheorem{dfn}[thm]{Definition}

\def\R{\hbox{$\mathbb R$}}

\def\Z{\hbox{$\mathbb Z$}}
\def\C{\hbox{$\mathbb C$}}

\def\ra{\rightarrow}
\def\xra#1{\xrightarrow{#1}}
\def \ie{\hbox{\it i.e.}}
\def\ainf{\mbox{$A_\infty$}}
\def\linf{\mbox{$L_\infty$}}
\def\tns{\otimes}
\def\ph{\varphi}

\def\la{\lambda}
\def\s#1{(-1)^{#1}}

\def\psa#1#2{\psi^{#1}_{#2}}

\def\zt{$\mathbb{Z}_2$}
\def\Cn#1#2{C^{#1}_{#2}}

\def\hom{\operatorname{Hom}}
\def\im{\operatorname{Im}}
\def\spn#1{\left\langle#1\right\rangle}

\def \De{D_{e_1}}
\def \Df{D_{f_s}}
\def\Da{D_a}
\def\Dx{D_x}
\def\thta{\la^{e_1f_s}-\la^{f_se_1}}
\def\Ch#1#2{\operatorname{Ch}^{#1}_{#2}}
\def\Hn#1#2{\operatorname{H}^{#1}_{#2}}
\def\He#1{\operatorname{H}^{#1}(E)}
\def\Hf#1{\operatorname{H}^{#1}(F)}
\def\larh{\la^{e_1}+\rho^{e_1}}

\begin{document}
\setlength{\multlinegap}{0pt}
\title[Nilpotent Algebra]{Computation of the cohomology of a family of \zt-graded Nilpotent Algebras}%

\address{}%
\email{}%

\thanks{}%
\date{\today}
\begin{abstract}
In this paper we compute the cohomology of certain special cases
of nilpotent algebras in a complex \zt-graded vector space of
arbitrary finite dimension. These algebras are generalizations of
the only two nontrivial complex 2-dimensional nilpotent algebras,
one from the moduli space of $1|1$-dimensional algebras, and the
other from the space of nongraded 2-dimensional algebras.
\end{abstract}
\maketitle
\section{Introduction}
The classification of complex associative algebras goes back at least to the
1870's, when Benjamin Peirce \cite{pie} gave a partial classification of
algebras of low dimension. At the time, deformation theory did not exist,
so Peirce used as a main tool the study of idempotents of the algebras. In
particular, he showed that a finite dimensional algebra which is not nilpotent
always contained an idempotent.

These ideas eventually led to the Fundamental
Structure Theorem for finite dimensional algebras, which says that a finite
dimensional algebra has a maximal nilpotent ideal, called its radical, and
the quotient of the algebra by its ideal is either trivial or semisimple.
The classification of simple algebras was given by Wedderburn's Theorem,
which says that a simple algebra is a tensor product of a matrix algebra
and a division algebra.
The division algebras over \C\ and \R\ are easily classified, so the semisimple
algebras can be easily determined. Therefore, to classify algebras using extensions,
it is important to understand the nilpotent algebras.

The infinitesimal deformations of associative algebras are classified by
Hochschild cohomology, which first appeared in \cite{hoch}. In \cite{gers},
the cohomology of associative algebras was studied,
and in \cite{gers1,gers2,gers3,gers4}, the deformation theory for associative
algebras was extensively studied. Gerstenhaber introduced the structure of
a graded Lie algebra on the Hochschild cochains, a rather mysterious fact whose
explanation was given in \cite{sta4}, where Stasheff realized that the space of
Hochschild cochains was the set of coderivations of the tensor coalgebra of the
underlying vector space of the associative algebra, which has a natural structure
of a graded Lie algebra.

In \cite{sta1,sta2}, Stasheff introduced the notion of an \ainf\ algebra, which
is a graded algebra, defined in terms of Hochschild cochains. The simplest way
to realize such an algebra is to consider a \zt-graded (or \Z-graded) vector
space, and then to construct the tensor coalgebra of its parity reversion (or
(de)suspension in the \Z-graded case). The coderivations of this tensor coalgebra
again form a graded Lie algebra, with respect to the internal grading induced only
by the underlying graded space, rather than the external grading which originally
was used in the Gerstenhaber construction. A codifferential is an odd (degree 1)
coderivation whose square is equal to zero, and Hochschild cohomology is induced
by the bracket (a fact noted by Gerstenhaber).

Using this approach, one of the authors developed a set of Maple worksheets which
allow the cohomology of \ainf\ algebras (and \linf\ algebras) to be computed
easily, at least in small degrees. However, the cohomology of an associative algebra,
even of low dimension, can be very complicated, but the computation of cohomology
given by a computer in small degrees often suggests a pattern which can be verified.
In fact, in \cite{bdhoppsw1,bdhoppsw2}, we gave a complete calculation of the cohomology
of every $1|1$ and ordinary 2-dimensional associative algebra.

As it turned out, there is exactly one nontrivial nilpotent algebra in the
$1|1$ case, and its cohomology pattern was complicated. Nevertheless, we were
eventually able to compute the cohomology of this nilpotent algebra. In the
case of ordinary 2-dimensional algebras, there is again exactly one nontrivial
nilpotent algebra, and while the pattern of its cohomology is simple, the verification
of  this pattern was as difficult as the $1|1$-dimensional case, and in fact,
the computation is parallel.

These two nilpotent algebras give rise to a family of nilpotent algebras in
dimension $r|s$. While studying $1|2$- and $2|1$-dimensional algebras, we noticed
that there was an interesting pattern for the cohomology of these nilpotent
algebras, and wondered if it would be possible to determine the complete
cohomology pattern. This paper contains that computation. There may be some
interesting applications of these algebras in the construction of some
\ainf\ extensions of these algebras.

The main result of this paper is the calculation of the dimensions of the cohomology
spaces $H^n$ for an $r|s$-dimensional space. Actually, we give a basis of the
cohomology space as well, which may be important in applications.

\section{Preliminaries}
If $V$ is a \zt-graded vector space, then an associative algebra structure
$m$ on $V$ is equivalent to an odd, quadratic codifferential
$d$ on the tensor coalgebra $T(W)$ of parity reversion $W=\Pi V$,
which is the space $V$ with the parity of elements reversed. In other words,
$d:W\tns W\ra W$ is a map which satisfies the following property,
\begin{equation*}
 d(d(a,b),c)+\s{a}d(a,d(b,c))=0
\end{equation*}

In this paper, $V$ is an $s|r$-dimensional vector space over $\C$,
so that its parity reversion $W$ is an $r|s$-dimensional.
The space of coderivations of $T(W)$ can be identified with the
space $C(W)=\hom(T(W),W)$ of cochains of $W$,
which is the direct product $C(W)=\prod C^n(W)$
of the spaces of $n$-cochains $C^n(W)=\hom(T^n(W),W)$.
Suppose that $W=\langle e_i,f_j\rangle_{i=1,j=1}^{r,s}$,
where the $e_i$ are even and the $f_j$ are odd.
For a multi-index $I=(i_1,\dots,i_n)$, where $i_k\in\{e_i,f_j\}$ and
$i\in\{e_i,f_j\}$, let $\ph^I_i:T^n(W)\ra W$ be given by
\begin{equation}
 \ph^I_i(J)=\delta^I_Ji, \nonumber
\end{equation}
where $J=j_1\tns\dots\tns j_n$ for the multi-index $J$,
and $\delta^I_J$ is the Kronecker delta symbol.
Then using Einstein summation notation,
we can represent any $n$-cochain $\ph$ in the form $\ph=\ph^I_ia^i_I$,
where $a^i_I\in\C$. In order to emphasize the
distinction between even and odd cochains,
we will often denote $\ph^I_i$ as $\psi^I_j$ when it is odd.

The space $C(W)$ is endowed with a natural Lie algebra structure.
Let $\ph^I_i\in C^n(W)$ and $\ph^J_j\in C^m(W)$, then
\begin{equation*}
 [\ph^I_i,\ph^J_j]=\ph^I_i\circ\ph^J_j-\s{(|I|+|i|)(|J|+|j|)}\ph^J_j\circ\ph^I_i,
\end{equation*}
where $|\cdot|$ represents the parity of the element.
Note that we typically drop the parity notation in a
power as it is clear what is meant.
Composition of elements is given by
\begin{equation*}
 \ph^I_i\circ\ph^J_j=
 \sum_{k=1}^{\ell(I)}\s{(i_1+\cdots+i_k)(J+j)}\delta^{i_k}_j\ph^{(I,J,k)}_i,
\end{equation*}
where $(I,J,k)$ is given by inserting $J$ into $I$ in
the place of the $k^{th}$ element of $I$;
\ie\ $(I,J,k)=(i_1,\cdots,i_{k-1},j_1,\cdots,j_{\ell(J)},i_{k+1},\cdots,i_{\ell(I)})$.

An odd coderivation $d$ is called a \emph{codifferential} if
\begin{equation*}
[d,d]=0.
\end{equation*}
We define the coboundary operator $D:C(W)\ra C(W)$ by $D(\ph)=[d,\ph]$.
Note that this is a differential on $C(W)$ since $D^2=0$.

Given a multi-index $J$, we define two maps $\la^J,\rho^J:C^n(W)\ra C^{n+\ell(J)}$
of parity $|J|$ and \emph{exterior degree} $\ell(J)$ by
\begin{align*}
\la^J\ph^I_i&=\ph^{JI}_i\\
\rho^J\ph^I_i&=\s{IJ}\ph^{IJ}_i.
\end{align*}

\section{Main Results}
\begin{thm}
 Let $W=\langle e_i,f_j\rangle_{i=1,j=1}^{r,s}$
 for $r,s\ge1$ and $d=\psa{e_1e_1}{f_s}$,
then $d$ is a codifferential and
\begin{equation*}
 h^n=\dim(H^n)=\left\{ \begin{array}{cc}
              (r+s-1)^{n+1} & n=0,1\mod4 \\
          (r+s-1)^{n+1}+1 & n=2,3\mod4
             \end{array}
         \right. \quad \text{for all }n
\end{equation*}
\label{thm-main}
\end{thm}
\begin{proof}
Clearly $d$ is a codifferential.
 \begin{align*}
D(\ph^I_{e_1})=&\ph^{e_1I}_{f_s}+\ph^{Ie_1}_{f_s}+(-1)^{I+1}\ph^I_{e_1}\psa{
e_1e_1}{f_s} &  \\
D(\ph^I_a)=&(-1)^{I+a+1}\ph^I_{a}\psa{e_1e_1}{f_s} &  a\ne e_1
 \end{align*}
Thus we see
\begin{align*}
 D=\De+\De^{f_s}:\Cn{n}{e_1}\ra&\Cn{n+1}{e_1}\oplus\Cn{n+1}{f_s} \\
 \Da:\Cn{n}{a}\ra&\Cn{n+1}{a} & a\ne e_1
\end{align*}
Since $D^2=0$ we have
\begin{align*}
 \De^2&=0\\
 \Df\De^{f_s}&=-\De^{f_s}\De
\end{align*}
We will first construct a space which will be shown to be exactly the
cohomology of the restricted $\Da$ operator. This space will then
be modified for $a=e_1$ and $a=f_s$.

Define the map $\theta=\thta$.
We claim that $\theta$, $\la^a$ and $\larh$ graded commute with the
restricted $D$ operator for $a\ne e_1,f_s$.
Since $\la^a$ and $\larh$ don't interact with the operator, they clearly commute.
For $\theta$ we compute
\begin{align*}
 \Da(\theta\ph^I_a)=&\Da((\thta)\ph^I_a)\\
  =&(-1)^{I+a+1}\ph^{e_1e_1e_1I}_a+
 (-1)^{I}\la^{e_1f_s}D(\ph^I_a)-(-1)^{I+a+1}\ph^{e_1e_1e_1I}_a\\
&-(-1)^{I}\la^{f_se_1}D(\ph^I_a)\\
=&-\theta D(\ph^I_a).
\end{align*}

Using these maps a space of non-trivial cocycles is constructed,
\begin{align*}
 \Ch{0}{x}=&\spn{\ph_x}\\
 \Ch{1}{x}=&\spn{\ph^{e_1}_x,\la^a\Ch{0}{x}|a\ne e_1,f_s}\\
 \vdots & \\
 \Ch{n}{x}=&\spn{\la^a\Ch{n-1}{x},\theta\Ch{n-2}x,
 (\larh)\la^a\Ch{n-2}x|a\ne e_1,f_s}
\end{align*}
Since the images of the maps $\la^a,(\larh)\la^a$
and $\theta$ are linearly independent,
$\Ch{n}{x}$ is a subspace of non-trivial cocycles.
Our goal is to show that this space is exactly
the cohomology of the restricted $D$ operator.

To accomplish this we claim that $\dim(\Ch{n}{x})=(r+s-1)^n$.
For $n=0$ this is trivial. Assume that we have shown for $k<n$, then
\begin{align*}
 \dim(\Ch{n}x)=&(r+s-2)(r+s-1)^{n-1}+(r+s-1)^{n-2}+(r+s-2)(r+s-1)^{n-2}\\
        =&(r+s-2)(r+s-1)^{n-1}+(r+s-1)^{n-1}\\
    =&(r+s-1)^n
\end{align*}

\begin{dfn}
 The spaces $B^n_a$ and $Z^n_a$ are the boundaries
 and cocycles of the $D_a$ operator,
 respectively. Also $b^n_a$ and $z^n_a$ are the relative
 dimensions of $B^n_a$ and $Z^n_a$.
\end{dfn}

We then have $H^n_a=B^n_a/Z^n_a$ and $h^n_a=dim(H^n_a)$, and $h^n_a=z^n_a-b^n_a$.
Also $z^n_a+b^{n+1}_a=(r+s)^n$.
When no confusion would result, we drop subscripts on $B^n_a$, $Z^n_a$ and $H^n_a$.

Since $\theta$ commutes with the restricted $D$ operator then we have,
\begin{equation*}
 \theta:B^n\ra B^{n+2}
\end{equation*}
We also observe that,
\begin{equation*}
 \Dx\la^a:\Cn{n}x\ra B^{n+2} \qquad a\ne e_1.
\end{equation*}
Since the images of these maps are disjoint we
see that $b^n\ge b^{n-2}+(r+s-1)(r+s)^n$.

We now claim that $b^{n+1}+b^n\ge (r+s)^{n+1}-(r+s-1)^{n+1}$.
By calculation we see that $b^0=0$ and $b^1\ge 1$, so this is true for $n=0$.
Assume we have shown for $k<n$, then
\begin{align*}
 b^{n+1}+b^n\ge& b^{n-1}+(r+s-1)(r+s)^n+b^n\\
            \ge& (r+s-1)(r+s)^n+(r+s)^{n-1}-(r+s-1)^{n-1}\\
        \ge& (r+s)^{n+1}-(r+s-1)^{n+1}
\end{align*}

Since $\Ch{n}{x}$ is a space of non-trivial cocycles we have
\begin{align*}
 (r+s-1)^n\le h^n=&z^n-b^n\\
                 =&(r+s)^n-b^{n+1}-b^n\\
         \le& (r+s)^n-((r+s)^n-(r+s-1)^n)\\
         =& (r+s-1)^n
\end{align*}
This shows that the cohomology associated with the restricted
$D$ operator is identifiable with $\Ch{n}x$.

We now extend this cohomology to the full $D$ operator.
Since the $D$ operator on the spaces $\Cn{n}{a}$ is the same as the
restricted $D$ operator for $a\ne e_1,f_s$, we are finished with these spaces.
However, the $\Cn{n}{e_1}$ and $\Cn{n}{f_s}$ spaces do interact.
To deal with this interaction we require additional terminology.

\begin{dfn}
 Define the a linear map $\tau:\Cn{n}{e_1}
 \oplus\Cn{n}{f_s}\ra\Cn{n}{e_1}\oplus\Cn{n}{f_s}$ by
 \begin{align}
  \tau(\ph^I_{e_1})&=\ph^I_{f_s} & \text{and}& &
  \tau(\ph^I_{f_s})&=\ph^I_{e_1} \nonumber
 \end{align}
\end{dfn}

\begin{dfn}
 We say $\ph\in\Ch{n}{e_1}$ \emph{extends} if there is an
 $\eta\in\Cn{n}{f_s}$ so that $D(\ph+\eta)=0$, moreover we say $\eta$ extends $\ph$.
 Also we say $\xi\in\Ch{n}{f_s}$ is trivial if there is an $\alpha\in\Cn{n}{e_1}$ so
 that $D(\alpha)=\xi$.
\end{dfn}

We now define the following spaces
\begin{align*}
 \Hn{n}{e}=&\frac{\ker(\De)}{\im(\De)} &
 \Hn{n}{f}=&\frac{\ker(\Df)}{\im(\Df)} &
 \Hn{n}{}=&\frac{\ker(D:\Cn{n}{e_1,f_s}
 \ra\Cn{n+1}{e_1,f_s})}{\im(D:\Cn{n-1}{e_1,f_s}\ra\Cn{n}{e_1,f_s})}.
\end{align*}
Note that $\Hn{n}{e}\cong\Ch{n}{e_1}$ and $\Hn{n}{f}\cong\Ch{n}{f_s}$.
We have we a well defined map $\Hn{n}{f}\ra\Hn{n}{}$, denote its image by $\Hf{n}$.
Define the space $\He{n}=\Hn{n}{}\setminus\Hf{n}$,
this space is isomorphic to the extendable subspace of $\Hn{n}{e}$.

Now we claim that the following sequence is exact,
\begin{equation}
 0\ra\He{n}\xra{\tau}\Hn{n}{f}\xra{\larh}\Hn{n+1}f\ra\Hf{n+1}\ra0.
 \label{eqn_seq}
\end{equation}

First note that $\De^{f_s}=(\larh)\tau$.
Let $\ph\in\He{n}$. Then $(\larh)\tau(\ph)=\De^{f_s}(\ph)=D(\ph)$ so the image of
$\ph$ in $\Hn{n+1}{f}$ is zero. So $\im(\tau)\subseteq\ker(\larh)$.
Suppose, on the other hand, that $\xi\in\ker(\larh)$ then $\De^{f_s}(\tau\xi)=(\larh)\xi=0$,
but this implies that $\tau\xi$ is extendable. Thus $\im(\tau)=\ker(\larh)$.
Let $\xi=(\larh)\eta$ for $\eta\in\Hn{n}{f}$,
but then $\xi=\De^{f_s}(\tau\eta)=D(\tau\eta)$.
Thus the image of $\xi$ in $\Hf{n+1}$ is zero, and
(\ref{eqn_seq}) is an exact sequence. Therefore, we have
\begin{equation}
 \dim(\Hf{n})=\dim(\He{n-1})+(r+s-2)(r+s-1)^{n-1}
 \label{eqn_dim}
\end{equation}

Next we claim that $\theta(\larh)\approx(\larh)\theta$ in $\Hn{n}{}$,
by ``$\approx$'' we mean they differ by a coboundary. Let $x\in\Hn{n}{f}$. Then
\begin{align*}
\theta(\larh)x\approx&(\la^{e_1f_se_1}-\la^{f_se_1e_1}+\rho^{e_1}
\theta)x\pm\Df(\la^{f_sf_s}x)\\
=&(\la^{e_1f_se_1}-\la^{f_se_1e_1}+\rho^{e_1}\theta+\la^{f_se_1e_1}-\la^{e_1e_1f_s})x\\
=&(\la^{e_1f_se_1}-\la^{e_1e_1f_s}+\rho^{e_1}\theta)x\\
\approx&(\larh)\theta x.
\end{align*}
Since $(\larh)\tau=\De^{f_s}$, we know that elements
in $\Ch{n}{f_s}$ which are in the image of $\larh$, are trivial.
Similarly, by the above, elements in the image of $\theta^m(\larh)$ are also trivial.

We encounter a problem when have an element of the form $\theta^n$.
However, this specific case has been studied extensively in \cite{bdhoppsw1}
and can be seen to be
\begin{center}
\begin{tabular}{|c|c|c|c|}
\hline
$n\mod 4$&$\dim(\He{n})$&$\dim(\Hf{n})$\\
\hline
0&0&1\\
1&1&0\\
2&1&1\\
3&1&1\\
\hline
\end{tabular}
\end{center}
This implies that, for certain values of $n$,
we subtract one from the dimensions of $\Hf{n}$ and $\He{n}$.

Using these restrictions we have the following bound on the dimension of $\Hf{n}$:
\begin{equation}
 \dim(\Hf{n})\le(r+s-1)^n-(r+s-2)\sum_{i=1}^{\left\lfloor\frac{n}2\right\rfloor}
 (r+s-1)^{n-2i}+\left\{\begin{array}{ll}
                        -1 & n=1\mod4 \\
                0 & \text{elsewise}
                       \end{array}\right.
 \label{eqn_cnt}
\end{equation}

Now we claim that if $x\in\Hn{n}{f}$ is trivial then $\tau(x)$ extends.
If $x$ is trivial then $x=aD(\alpha)+bD(\beta)$, where $a,b\in\C$,
$\alpha\in\Cn{n-1}{e_1}$ and $\beta\in\Cn{n-1}{f_s}$, note that $a\ne0$. Then
\begin{align*}
 \De^{f_s}(\tau(x))=&(\larh)\tau\tau(x)\\
             =&(\larh)(x)
\end{align*}
Which is a coboundary since $\larh :B^n_{f_s}\ra B^{n+1}_{f_s}$.
Thus $\De^{f_s}(\tau(x))$ trivial and $\tau(x)$ extends.

An immediate consequence of this and the that occasionally $\theta^n$
will both extend and be non-trivial, we have
\begin{equation}
 \dim(\He{n})+\dim(\Hf{n})\ge(r+s-1)^n+\left\{\begin{array}{ll}
                                               0 & n=0,1\mod4 \\
                           1 & \text{elsewise}
                                              \end{array}\right.
 \label{eqn_rel}
\end{equation}

Combining equations (\ref{eqn_dim}) and (\ref{eqn_rel}) we arrive at
\begin{equation*}
 \dim(\Hf{n+1})+\dim(\Hf{n})\ge(r+s-1)^{n+1}+\left\{\begin{array}{ll}
                                               0 & n=0,1\mod4 \\
                           1 & \text{elsewise.}
                                              \end{array}\right.
\end{equation*}
Using this and (\ref{eqn_cnt}) we discover that this is actually and
equality which we restate as
\begin{equation}
  \dim(\He{n})+\dim(\Hf{n})=(r+s-1)^n+\left\{\begin{array}{ll}
                                               0 & n=0,1\mod4 \\
                           1 & \text{elsewise.}
                                              \end{array}\right.
\end{equation}

Finally, to obtain the total dimension of the cohomology we sum the individual pieces,
\begin{align}
 h^n &= \sum_{x\in\{e_i,f_j\}} h^n_x \nonumber \\
     &= (r+s-2)(r+s-1)^n + (r+s-1)^n+\left\{\begin{array}{ll}
                                               0 & n=0,1\mod4 \\
                           1 & \text{elsewise}
                                              \end{array}\right. \nonumber \\
     &= (r+s-1)^{n+1}+\left\{\begin{array}{ll}
                          0 & n=0,1\mod4 \\
               1 & \text{elsewise}
                         \end{array}\right.
\end{align}

\end{proof}

\begin{thm}
 Let $W=\langle e_i,f_j\rangle_{i=1,j=1}^{r,s}$ for $s\ge2$ and $d=\psa{f_1f_1}{f_s}$,
then $d$ is a codifferential and
\begin{equation*}
 h^n=(r+s-1)^{n+1}+1 \qquad \text{for all }n
\end{equation*}
\label{thm-main2}
\end{thm}

\begin{proof}
 This proof is analogous to the proof of Theorem
 \ref{thm-main} with the following exceptions.
 Define $\theta = \la^{f_1f_s}+\la^{f_sf_1}$ and
 \begin{align*}
 \Ch{0}{x}=&\spn{\ph_x}\\
 \Ch{1}{x}=&\spn{\ph^{f_1}_x,\la^a\Ch{0}{x}|a\ne f_1,f_s}\\
 \vdots & \\
 \Ch{n}{x}=&\spn{\la^a\Ch{n-1}{x},\theta\Ch{n-2}x,(\la^{f_1}+\rho^{f_1})\la^a\Ch{n-2}x|a\ne f_1,f_s}
\end{align*}
then, by \cite{bdhoppsw2}, we have
 \begin{center}
\begin{tabular}{|c|c|c|c|}
\hline
$n\mod 4$&$\dim(\He{n})$&$\dim(\Hf{n})$\\
\hline
0&1&1\\
1&1&1\\
2&1&1\\
3&1&1\\
\hline
\end{tabular}
\end{center}

Thus the cohomology will be exactly as claimed.

\end{proof}

\bibliographystyle{amsplain}

\providecommand{\bysame}{\leavevmode\hbox
to3em{\hrulefill}\thinspace}
\providecommand{\MR}{\relax\ifhmode\unskip\space\fi MR }
\providecommand{\MRhref}[2]{%
  \href{http://www.ams.org/mathscinet-getitem?mr=#1}{#2}
} \providecommand{\href}[2]{#2}

\end{document}